\documentclass{article}
\usepackage{amssymb,amsmath}
\usepackage{graphicx}
\usepackage{ dsfont }

\def\qed{\hfill {\hbox{${\vcenter{\vbox{               
   \hrule height 0.4pt\hbox{\vrule width 0.4pt height 6pt
   \kern5pt\vrule width 0.4pt}\hrule height 0.4pt}}}$}}}

\def\tr{\triangleright}

\def\hat{\widehat}

\newtheorem{definition}{Definition}

\newtheorem{example}{Example}

{\qed\par\medskip}

\date{}

\title{\Large \textbf{Quotient Quandles and the Fundamental Latin Alexander 
Quandle}}

\author{
Sam Nelson \footnote{knots@esotericka.org} \and 
Sherilyn Tamagawa\footnote{stamagaw2510@ScrippsCollege.edu}}

\begin{document}
\maketitle

\begin{abstract}
Defined in \cite{J,M}, the fundamental quandle is a complete invariant of 
oriented classical knots. We consider invariants of knots defined from
quotients of the fundamental quandle. In particular, we introduce the
fundamental Latin Alexander quandle of a knot and consider its Gr\"obner
basis-valued invariants, which generalize the Alexander polynomial. We show
via example that the invariant is not determined by the generalized 
Alexander polynomial for virtual knots.
\end{abstract}

\parbox{5.5in} {\textsc{Keywords:} quandles, Alexander quandles, Latin 
quandles, knot invariants, Gr\"obner bases
\smallskip

}

\section{\large \textbf{Introduction}}

In \cite{J,M} Joyce and Matveev introduced the algebraic structures known
as \textit{quandles} or \textit{distributive groupoids}. In particular,
the \textit{fundamental quandle} of an oriented classical knot was shown to
determine the knot group and the peripheral subgroup and thus the knot complement
up to homeomorphism, yielding a complete invariant of oriented classical
knots. In \cite{J,W} quotients of the fundamental quandle, including the
fundamental involutory quandle and the fundamental involutory abelian quandle
were studied, including some connections to the Alexander invariant. 
In particular, Joyce showed that the fundamental involutory abelian quandle 
of a knot is always finite with cardinality equal to the determinant of the 
knot, while Winker showed that some knots have infinite involutory 
(non-abelian) quandle.

In this paper we consider some quotients of the fundamental quandles of 
classical and virtual knots and describe an algorithm which can sometimes 
reveal when a quotient of the fundamental quandle of a knot is finite.
Showing that a given quotient is infinite is harder for general quotient 
quandles, but is simpler for quotients of the fundamental Alexander quandle 
of a knot, which has a module structure. We introduce the \textit{fundamental 
Latin Alexander quandle} of a knot, a generalization of the Alexander quandle 
with coefficients in an extension ring such that the resulting quandle is 
Latin. From this new structure we define Gr\"obner basis-valued invariants 
akin to those defined in \cite{CHN}. We include an example which shows that 
the new invariant is not determined by the generalized Alexander polynomial 
for virtual knots.

The paper is organized as follows. In Section \ref{Q} we review the basics
of quandles. In Section \ref{QFQ} we consider some quotients of the fundamental
quandle. In Section \ref{FLAQ} we define the Fundamental Latin Alexander quandle
and the Fundamental Latin Alexander Gr\"obner (FLAG) invariants, including 
computations of the $\mathrm{FLAG}_1$ invariant for all classical knots with
up to eight crossings. We end in Section \ref{QFR} with some questions for
future research.

\section{\large \textbf{Quandles}}\label{Q}

We begin with a definition (see \cite{J,M,FR}).

\begin{definition}\textup{
A \textit{quandle} is a set $Q$ with an operation $\tr:Q\times Q\to Q$ 
satisfying for all $x,y,z\in Q$
\begin{itemize}
\item[(i)] $x\tr x=x$,
\item[(ii)] The map $f_y:Q\to Q$ defined by $f_y(x)=x\tr y$ is a bijection, and
\item[(iii)] $(x\tr y)\tr z=(x\tr z)\tr(y\tr z)$.
\end{itemize}
The inverse of $f_y$ defines another operation called the \textit{dual quandle}
operation $f_y^{-1}(x)=x\tr^{-1} y$.
}\end{definition}

It is a straightforward exercise to show that $Q$ forms a quandle under the 
dual quandle operation and that the two operations mutually right-distribute, 
i.e., we have
\begin{eqnarray*}
(x\tr y)\tr^{-1} z & = & (x\tr^{-1}z)\tr(y\tr^{-1}z)\\
(x\tr^{-1} y)\tr z & = & (x\tr z)\tr^{-1}(y\tr z).
\end{eqnarray*}

\begin{example}\textup{
Any $\mathbb{Z}$-module $A$ is a quandle under the operation
\[x\tr y=2y-x.\]
In particular, the dual quandle operation is the same as the original quandle
operation, i.e. $x\tr^{-1} y=x\tr y$. Quandles with this property are called
\textit{involutory} since the maps $f_y$ are involutions.
}\end{example}

\begin{example}\textup{
Let $G$ be any group. Then $G$ is a quandle under $n$-fold conjugation
\[x\tr y=y^{-n}xy^n\]
and under the \textit{core} operation
\[x\tr y=yx^{-1}y.\]
The set $G$ with these quandle structures is called $\mathrm{Conj}_n(G)$ and
$\mathrm{Core}(G)$ respectively.
}\end{example}

\begin{example}\textup{Any module $M$ over the ring 
$\Lambda=\mathbb{Z}[t^{\pm 1}]$ is a quandle under the operation 
\[x\tr y=tx+(1-t)y\]
called an \textit{Alexander quandle}. More generally, if $A$ is any
abelian group and $t\in\mathrm{Aut}(A)$ is an automorphism of abelian groups,
then $A$ is an Alexander quandle under the operation above where $1$ is the
identity map.
}\end{example}

\begin{example}\textup{
Let $K$ be a link in $S^3$ and $N(K)$ a regular neighborhood of $K$. Then the
\textit{fundamental quandle} of $K$ is the set of homotopy classes of paths 
in $S^3\setminus N(K)$ from a base point to $N(K)$ such that the initial point
stays fixed at the base point while the terminal point is free to wander on 
$N(K)$. The quandle operation is then given by setting $x\tr y$ to the homotopy
class of the path given by first following $y$, then going around a canonical
meridian on $N(K)$ linking $K$ once, then going backward along $y$, then 
following $x$ as illustrated. See \cite{J} for more.
\[\includegraphics{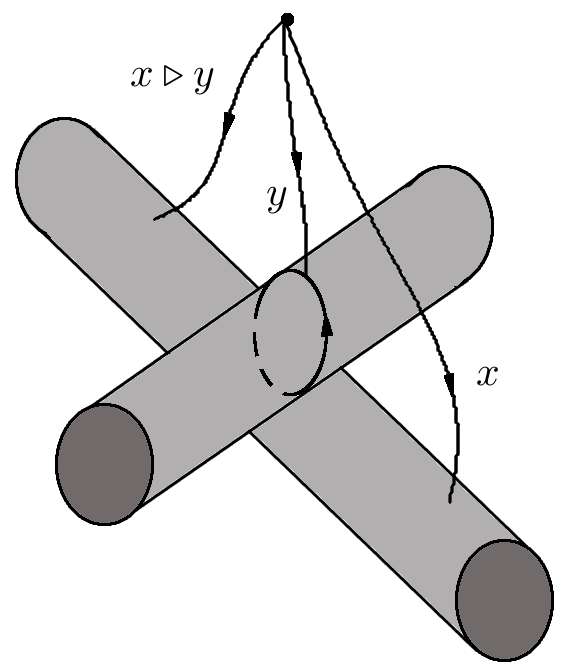}\]
}\end{example}

\begin{example}\label{ex:f8}
\textup{The knot quandle can also be expressed combinatorially
with a presentation by generators and relations
as the set of equivalence classes of quandle words in a set of generators 
corresponding to arcs in a diagram of $K$ under the equivalence relation
generated by the quandle axioms together with the \textit{crossing relations}
\[\includegraphics{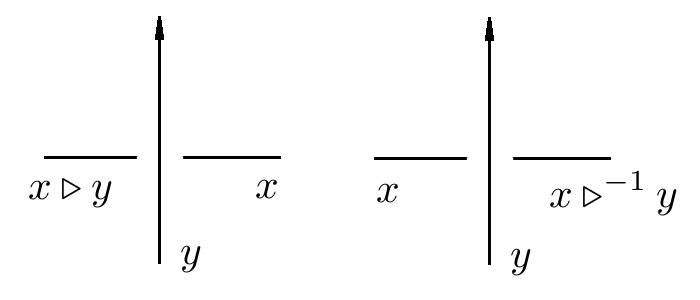}\]
For example, the figure 8 knot $4_1$ below has the listed quandle presentation
\[\includegraphics{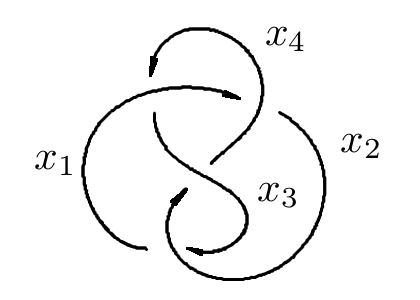} \quad\raisebox{0.5in}{$\langle x_1,x_2,x_3,x_4\ 
|\ x_2\tr x_4=x_1, x_3\tr x_2=x_4, x_4\tr x_2=x_1, x_2\tr x_3=x_4\rangle.$}\]
}\end{example}

Say a relation in a quandle presentation is \textit{short} if it has the
form $x_i\tr x_j=x_k$ for $x_i,x_j,x_k$ generators. Then we observe that 
every finitely presented quandle $Q$ has a presentation in which every 
relation is short, since we can add new generators $x_k$ and short form 
relations abbreviating subwords of the form $x_i\tr x_j$ to $x_k$
as needed until all relations are short. If our generators are numbered 
$\{x_1,\dots,x_n\}$, then we can express a short form presentation with a 
matrix whose row $i$ column $j$ entry is $k$ if $x_i\tr x_j=x_k$ and $0$ 
otherwise; we will call this a \textit{presentation matrix} for $Q$. If a 
presentation matrix for $Q$ has no zeros, then it expresses the complete 
operation table for $Q$, and $Q$ is a finite quandle. 

\begin{example}\textup{
The figure 8 knot in example \ref{ex:f8} above has presentation matrix
\[\left[\begin{array}{rrrr}
0 & 0 & 0 & 0 \\
0 & 0 & 4 & 1 \\
4 & 1 & 0 & 0 \\
0 & 0 & 0 & 0 \\
\end{array}\right].\]
}\end{example}

\section{\textbf{Quotients of the Fundamental Quandle}}\label{QFQ}

Knot quandles are generally infinite. However, it is observed in \cite{J} and
later in \cite{W} that the \textit{involutory} version of the fundamental
quandle of a knot is often finite, and the fundamental \textit{abelian 
involutory quandle} of a knot is always finite with order equal to the 
determinant of the knot, i.e. the absolute value of the Alexander polynomial 
evaluated at $-1$. The fundamental involutory quandle of a knot can be 
understood as the result of adding a fourth axiom which says 
\begin{itemize}
\item[(iv)] $x\tr y=x\tr^{-1}y$ for all $x,y\in Q$
\end{itemize}
or equivalently, replacing the second quandle axiom with
\begin{itemize}
\item[(ii)'] $(x\tr y)\tr y=x$ for all $x,y\in Q$;
\end{itemize}
the fundamental abelian involutory quandle is then obtained by adding an 
additional axiom which says
\begin{itemize}
\item[(v)] $(x\tr y)\tr(z\tr w)=(x\tr z)\tr(y\tr w)$ for all $x,y,z,w\in Q$.
\end{itemize}

We can verify that the fundamental involutory quandle of the figure eight knot
is finite with cardinality 5 by observing that moves of the following types do
not change the quandle presented by a presentation matrix:
\begin{itemize}
\item[(i)] Filling in a zero with a value obtained as a consequence of the 
axioms and other relations,
\item[(ii)] Filling in a zero with a number defining a new generator and adding
a row and column of zeroes corresponding to the new generator
\item[(iii)] Deleting a row and column and replacing all instances of the 
larger generator with the smaller one when two generators are found to be equal,
taking care to note any new equalities of generators implied.
\end{itemize}

This gives us an algorithm for filling in the complete operation table of a
finitely presented quandle: first, fill in all zeroes determined by
consequences of the axioms and keep a list of any pairs of equal generators,
reducing the presentation by eliminating redundant generators when possible. 
Next, if any zeroes remain, choose one to assign to a new generator and repeat
the process. This procedure may or may not terminate -- if the presented 
quandle is infinite, the process can never terminate, but even if the quandle 
finite then the speed of termination depends a great deal on the choice of 
zeroes for replacement. On the other hand, when the process does terminate, the
result is a sequence of Tietze moves showing that the presented quandle is 
finite.

\begin{example}\label{ex:f82}\textup{
Let us use the above procedure to verify that the figure eight knot has 
fundamental involutory quandle of cardinality 5. We start with the presentation
matrix from example \ref{ex:f8} and fill in the zeroes as determined by the
involutory quandle axioms:
\[\left[\begin{array}{rrrr}
0 & 0 & 0 & 0 \\
0 & 0 & 4 & 1 \\
4 & 1 & 0 & 0 \\
0 & 0 & 0 & 0 \\
\end{array}\right]
\longrightarrow
\left[\begin{array}{rrrr}
1 & 3 & 0 & 2 \\
0 & 2 & 4 & 1 \\
4 & 1 & 3 & 0 \\
3 & 0 & 2 & 4 \\
\end{array}\right].
\]
For example, quandle axiom (i) says $x_i\tr x_i=x_i$, so the diagonal elements
are filled in with their row numbers; the involutory condition says 
that since $x_3\tr x_2=x_1,$ we have
$x_1\tr x_2=(x_3\tr x_2)\tr x_3=x_3$, etc. Note that we still have some zeroes
which cannot be filled in from the axioms; thus, we need to choose a zero 
to assign a new generator $x_5$ -- say we set $x_5=x_1\tr x_3$. Then we have 
presentation matrix below which completes via the involutory quandle axioms to
the matrix on the right
\[
\left[\begin{array}{rrrrr}
1 & 3 & 5 & 2 & 0\\
0 & 2 & 4 & 1 & 0\\
4 & 1 & 3 & 0 & 0\\
3 & 0 & 2 & 4 & 0\\
0 & 0 & 0 & 0 & 5
\end{array}\right]
\longrightarrow
\left[\begin{array}{rrrrr}
1 & 3 & 5 & 2 & 4\\
5 & 2 & 4 & 1 & 3\\
4 & 1 & 3 & 5 & 2\\
3 & 5 & 2 & 4 & 1\\
2 & 4 & 1 & 3 & 5
\end{array}\right]
\]
and the fundamental involutory quandle of $4_1$ has 5 elements.
}\end{example}

The involutory and abelian conditions have geometric motivations: the 
involutory condition comes from considering unoriented knots, while the abelian
condition is the condition required for the set of quandle homomorphisms 
from the knot quandle to $Q$ to inherit a natural quandle structure 
(see \cite{CN} for more). Nevertheless, we can consider these quotient
quandles to be simply the result of imposing algebraic conditions on the
fundamental quandle of a knot. Any such choice of conditions results in 
a quandle-valued knot invariant, and for each such invariant we can ask 
whether the resulting quandle is finite. In \cite{J,W} the generalizations
of the involutory condition to higher numbers of operations, e.g.
\[(\dots((x\tr y)\tr y)\dots) \tr y=x\]
were considered, with the notable result that the square knot and granny knot
have non-isomorphic 4-quandles (i.e., quotients in which we set 
$(((x\tr y)\tr y)\tr y)\tr y=x$), despite having isomorphic knot groups.

We considered several examples of algebraic axioms and used our procedure 
outlined above, implemented in \texttt{Python}, to search for examples of
knots whose fundamental quandles had finite quotients when the axioms were 
imposed. These included:
\begin{itemize}
\item \textit{Anti-abelian axiom}: $(x\tr y)\tr(z \tr w)=(w\tr y)\tr (z\tr x)$,
\item \textit{Left distributive axiom}: $x\tr(y\tr z)=(x\tr y)\tr (x\tr z)$
\item \textit{Commutative operator axiom}: $x\tr(y\tr z)=x\tr(z\tr y)$
\item \textit{Latin axiom}: $x\tr y=x\tr z \Rightarrow y=z$
\end{itemize}
both in combination with the abelian and involutory axioms and alone.
Some combinations are redundant; for instance, the abelian condition implies 
left distributivity. Curiously, we found that many of the above conditions 
yield the same results, with most knots of small crossing number having either
trivial one-element quotient quandles or the three-element quandle structure
$\mathbb{Z}_3$ with $x\tr y= 2x-y$.


\begin{example}\textup{
Of the classical knots with seven or fewer crossings, $3_1, 6_1, 7_4$,  
and $7_7$ have anti-abelian involutory quandles with three elements
$$\left[ \begin{array}{rrr}
1 & 3 & 2 \\
3 & 2 & 1 \\
2 & 1 & 3
\end{array}\right],$$
while the rest have the trivial one-element quandle.}
\end{example}

Recall that a \textit{virtual knot} is an equivalence class of oriented Gauss 
codes under the equivalence relation determined by the Gauss code Reidemeister
moves.  It is standard practice to draw
virtual knots with extra \textit{virtual crossings}, circled self-intersections
representing non-planarity of Gauss codes; these virtual crossings interact with 
classical crossings via the \textit{detour move}, 
\[\includegraphics{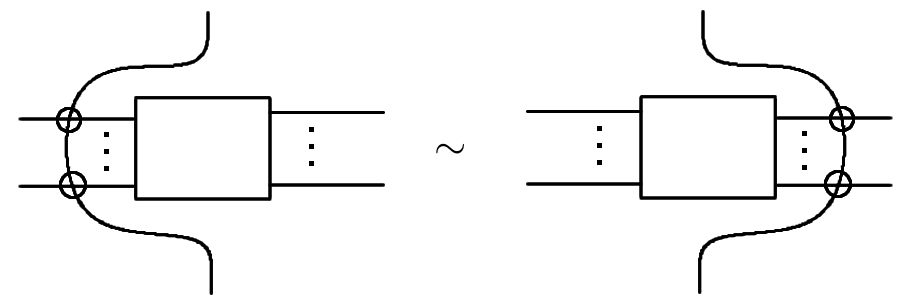}\]
which says we can redraw any 
arc with only virtual crossings in its interior as any other arc with only virtual 
crossings in its interior. Virtual knots may be understood as equivalence 
classes of knots in thickened oriented surfaces $\Sigma\times [0,1]$ modulo
stabilization. See \cite{K,KM} for more about virtual knots.

\begin{definition}\textup{
Let $Q$ be a quandle and $v:Q\to Q$ a bijective map.
We say that $(Q,v)$ is a \textit{virtual quandle} if $v$ satisfies
\[v(x\tr y)=v(x)\tr v(y),\]
i.e., a virtual quandle is a quandle with a choice of automorphism.
If $Q$ is an involutory quandle, then $(Q,v)$ is an \textit{involutory virtual 
quandle} if $(Q,v)$ is a virtual quandle and $v$ is an involution, i.e. if
$v(v(x))=x$ for all $x\in Q$.
}\end{definition}

Let $K$ be a virtual knot. The \textit{fundamental involutory virtual quandle}
of $K$ is the virtual quandle with presentation consisting of one generator
for each portion of $K$ containing only overcrossings (that is, we divide
$K$ at classical undercrossings and at virtual crossings) with relations as
pictured together with the involutory quandle axioms.
\[\includegraphics{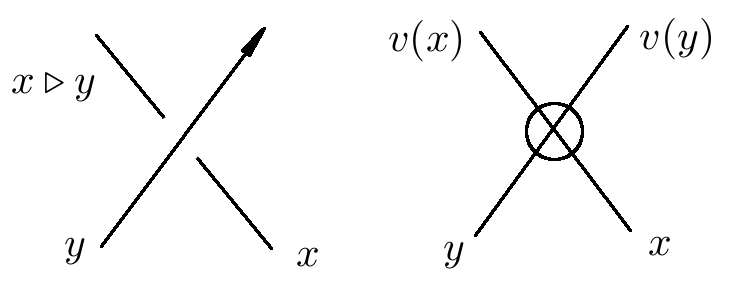}\]

\begin{example}\textup{
The virtual knots $3.2, 3.3, 3.4, 3.5, 4.2, 4.4, 4.5, 4.9,$ and $4.{43}$ all 
have anti-abelian involutory virtual quandle
\[\left[ \begin{array}{rrr|r}
1&3&2&2\\
3&2&1&1\\
2 &1&3&3 
\end{array}\right]\]
and the virtual knot $4_3$ has the trivial involutory virtual quandle.
}\end{example}

\begin{example}\textup{
However, while many knots seem to have either the trivial one-element 
anti-abelian involutory virtual quandle or the three-element 
anti-abelian involutory virtual quandle structure above, not all of 
them do.  For example, we found via \texttt{Python} computations that 
the anti-abelian quandle of  the virtual knot $3_7$ is}
\[\includegraphics{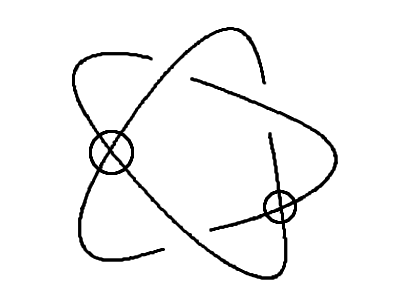}\]

\noindent\textup{
\scalebox{.74}{
$\left[ \begin{array}{rrrrrrrrrrrrrrrrrrrrrrrrrrr|r}
1& 7& 6& 8& 9& 3& 2& 4& 5& 11& 10& 13& 12& 15& 14&17& 16& 20& 21 &1& 19& 23& 22& 25& 24& 27& 26&2\\
 7&2& 12& 5&4&14&1& 10& 11&8&9& 3&15&6&13& 18& 21&16&20& 19 &17& 24& 26& 22& 27& 23& 25&1\\
 6&12&3& 16& 22&1&15& 26&21&20& 25& 2&14& 13&7& 4&27& 24&23& 10& 9& 5& 19& 18& 11& 8& 17&4\\
 8&5& 16&4& 2& 27& 11& 1& 10& 9& 7&24& 21& 20&23& 3& 26& 22 &25& 14& 13& 18& 15&12& 19& 17& 6&3\\
 9& 4 &22& 2& 5& 19& 10& 11& 1 &7& 8& 18& 26& 25& 17&24&15 &12& 6& 27&23 &3& 21& 16& 14 &13& 20&6 \\
3& 14& 1& 27 &19& 6& 13& 17& 23& 24 &18& 15& 7& 2& 12& 26& 8& 11& 5&25& 22& 21& 9& 10& 20& 16& 4&5 \\
2& 1& 15& 11& 10& 13&7&9& 8& 5& 4&14& 6& 12& 3 &19&20& 21& 16& 17& 18 &27 &25&26& 23& 24& 22&7\\
 4 &10& 26& 1&11&17&9&8& 7& 2&5& 19& 25&22& 18& 27& 6&15&12& 23& 24& 14 &20& 21&13& 3& 16&12\\
 5&11&21&10&1 &23& 8& 7& 9& 4 &2&27& 20&16& 24 &14&25& 26& 22&13& 3&19& 6&15&17&18& 12&14 \\
11& 8& 20&9& 7&24&5& 2 &4& 10& 1& 23& 16&21& 27 &13& 22& 25& 26& 3& 14& 17&12& 6& 18 &19& 15&13\\
 10&9& 25& 7& 8& 18& 4& 5& 2& 1& 11& 17& 22& 26& 19& 23 &12& 6& 15& 24& 27 &13 &16&20& 3 &14& 21&15\\
 13&3& 2 &24& 18& 15& 14& 19& 27& 23& 17 &12& 1& 7& 6& 22& 11& 5& 8& 26& 25& 16& 10& 4 &21& 20& 9&8\\
12& 15& 14& 21& 26& 7& 6& 25& 20 &16& 22& 1 &13& 3& 2& 10& 23& 27& 24& 9& 4&11 &17& 19& 8& 5& 18&10\\
 15&6& 13& 20&25& 2&12&22&16&21& 26& 7&3&14&1& 9& 24& 23& 27& 4& 10& 8& 18&17& 5& 11& 19&9 \\
14 &13& 7& 23& 17&12& 3& 18& 24&27&19& 6& 2&1&15&25&5&8&11&22& 26&20&4&9&16&21& 10&11\\
 17 &18& 4& 3&24& 26&19& 27 &14&13& 23& 22& 10& 9& 25& 16& 1& 2& 7&21&20&12&11& 5&15&6& 8&16\\
 16 &21&27& 26& 15&8&20& 6&25&22&12&11& 23& 24&5& 1&17& 19& 18&7&2& 10&13& 14& 9& 4& 3&18\\
20& 16&24&22&12&11&21&15& 26&25&6& 5 &27&23& 8& 2 &19 &18& 17& 1& 7& 4& 14& 3 &10& 9& 13& 17\\
21& 20& 23& 25& 6& 5& 16&12& 22& 26&15 &8& 24&27& 11&7& 18&17& 19& 2&1&9&3&13& 4&10& 14&19\\
 18 &19&10& 14& 27 &25& 17&23&13&3& 24 &26& 9& 4&22& 21& 7& 1 &2& 20& 16& 15& 8& 11& 6&12& 5&21\\
 19&17& 9& 13 &23& 22&18& 24& 3& 14&27&25& 4& 10 &26 &20& 2& 7&1& 16& 21&6& 5&8& 12& 15&11 &20\\
23& 24& 5& 18&3& 21& 27&14&19& 17& 13& 16& 11& 8& 20& 12& 10& 4 &9& 15& 6&22& 1& 2&26&25& 7&27\\
 22&26 &19&15& 21& 9& 25& 20& 6& 12&16& 10&17&18& 4& 11& 13& 14&3&8& 5& 1& 23& 27&7&2& 24&25\\
 25&22&18& 12&16&10& 26&21&15& 6&20& 4& 19& 17& 9& 5& 14&3& 13& 11& 8 &2& 27&24& 1& 7& 23&26\\
 24& 27&11&19&14& 20& 23&13&17&18&3& 21&8&5& 16& 15& 9& 10&4 &6& 12&26&7&1&25& 22&2&23 \\
27&23&8&17&13&16 &24&3& 18&19&14&20& 5&11& 21& 6& 4&9&10& 12 &15& 25&2&7&22& 26& 1&24\\
26& 25&17&6&20&4& 22&16&12 &15& 21& 9& 18&19&10&8&3&13& 14&5&11 &7& 24&23&2&1& 27&22
\end{array} \right]$}.}
\end{example}

\section{\large \textbf{Fundamental Latin Alexander Gr\"obner Invariants}}\label{FLAQ}

Let $K$ be a knot or link. The \textit{fundamental Alexander quandle} $FAQ(K)$
is the $\Lambda$-module generated by generators corresponding to arcs in a
diagram of $K$ with Alexander quandle operations at the crossings. As a 
$\Lambda$-module, the fundamental Alexander quandle of a knot is the classical
Alexander invariant.

Let $R$ be a polynomial ring and $M$ an $R$-module with presentation matrix 
$P\in M_{m,n}(R)$, i.e. the rows of $P$ correspond to generators of $M$ and 
the rows of $P$ express relations defining $M$. The \textit{$k$th elementary 
ideal} $I_k$ of $M$ is the ideal in $R$ generated by the $(n-k)$ (or $m-k$ if 
$m>n$) minors of $P$. It is a standard result (see \cite{L} for instance) that
changes to $P$ reflecting Tietze moves in the presentation of $M$ do not
change the elementary ideals, and hence these ideals are invariants of $M$.

\begin{example}\textup{Let $K$ be a knot and $P$ a $\Lambda$-module 
presentation matrix of
the Alexander quandle of $K$. The \textit{$k$th Alexander Polynomial} of $K$ 
is any generator $\Delta_k$ of the smallest principal ideal of $\Lambda$ 
containing the $k$th elementary ideal of $P$. Note that
$\Delta_k$ is defined only up to multiplication by units in $\Lambda$. In
particular, $\Delta_0=1$ for classical knots $K$, and $\Delta_1$ is often
called the Alexander polynomial.
}\end{example}

Recall that a quandle is \textit{Latin} if in addition to the 
right-invertibility required by quandle axioms 2, we also have 
left-invertibility. That is, a quandle $Q$ is Latin if it satisfies the axiom
\begin{itemize}
\item[(iv)] For every $x,y\in Q$, there is a unique $z\in Q$ such that 
$x\tr z=y$
\end{itemize}
or equivalently
\begin{itemize}
\item[(iv')] For every $x\in Q$, the map $f_x:Q\to Q$ defined by $f_x(y)=x\tr y$
is bijective.
\end{itemize}
A finite quandle is Latin if and only its operation table forms a Latin square,
i.e. if every row and column is a permutation of the elements of $Q$.

\begin{example}\textup{
An Alexander quandle is Latin iff $1-t$ is invertible. For instance, the 
Alexander quandle structure on $\mathbb{Z}_3$ with 
$t\in\mathrm{Aut}(\mathbb{Z}_3)$ given by multiplication by $2$ is Latin, 
while the Alexander quandle structure on $\mathbb{Z}_4$ with 
$t\in\mathrm{Aut}(\mathbb{Z}_4)$ given by multiplication by $3$ is not Latin:
\[\begin{array}{r|rrr}
\tr & 0 & 1 & 2  \\\hline
0 & 0 & 2 & 1 \\
1 & 2 & 1 & 0 \\
2 & 1 & 0 & 2
\end{array}
\quad\quad
\begin{array}{r|rrrr}
\tr & 0 & 1 & 2 & 3  \\\hline
0 & 0 & 2 & 0 & 2 \\
1 & 3 & 1 & 3 & 1 \\
2 & 2 & 0 & 2 & 0 \\
3 & 1 & 3 & 1 & 3
\end{array}.\]
}\end{example}

The element $(1-t)\in\Lambda$ is not invertible in $\Lambda$; its ``natural''
inverse is the Laurent series $1+t+t^2+\dots$. We prefer to stick to polynomial
rings, so we define the Fundamental Latin Alexander Quandle of an oriented 
link $L$ in the following way: let $\hat\Lambda=\mathbb{Z}[t,t^{-1}s,s^{-1}]$
where the variables $t^{-1}$ and $s^{-1}$ are new formal variables, not (yet)
inverses for $t$ and $s$, and then define the quotient ring 
$\hat\Lambda'=\hat\Lambda/(ss^{-1}-1,tt^{-1}-1,1-t-s)$. Then we define the
\textit{Fundamental Latin Alexander Quandle} of $L$, $FLAQ(L)$, to be the
$\hat\Lambda'$-module generated by a set of generators corresponding to
arcs in a diagram of $L$ with relations of the form $w=tx+sy$
at crossings as depicted below. Equivalently,  $FLAQ(L)$ can be regarded 
as the Alexander quandle of the knot with coefficients in the extension ring 
$\Lambda[(1-t)^{-1}]$ of $\Lambda=\mathbb{Z}[t^{\pm 1}]$.

In \cite{CHN}, Gr\"obner basis-valued invariants of knots and link were
defined from the Alexander biquandle by considering the pullback ideals of the
elementary ideals to a standard (non-Laurent) polynomial ring, then taking
the Gr\"obner basis of the resulting ideal with respect to a choice of 
monomial ordering. We can apply the same idea here to get a new Gr\"obner 
basis-valued invariant which we call the \textit{Fundamental Latin Alexander 
Gr\"obner} invariant, denoted $FLAG(L)$.

\begin{definition}\textup{
Let $L$ be an oriented link, $\hat\Lambda=\mathbb{Z}[t,t^{-1},s,s^{-1}]$ a 
four variable polynomial ring, and $P$ the coefficient matrix of the 
homogeneous system of linear equations with variables corresponding to arcs
in a diagram of $L$ and equations at crossings as depicted.
\[\includegraphics{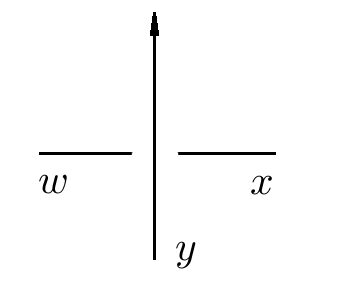}\quad\raisebox{0.5in}{$tx+sy=w$}\]
Then the $k$th FLAG ideal of $L$ is the ideal $I_k$ in 
$\hat\Lambda$ generated by the generators of the $k$th elementary ideal of $P$ 
and the polynomials $ss^{-1}-1,tt^{-1}-1$ and $1-t-s$. Given a choice of 
monomial ordering $<$ of the variables $s,s^{-1},t,t^{-1}$, the resulting 
Gr\"obner basis of $I_k$ is the $k$th FLAG invariant of $L$, denoted 
$FLAG^{<}_k(L)$.}
\end{definition}

The FLAG invariant contains in general more information that the Alexander 
polynomial; for instance, the number of elements of the FLAG basis with respect 
to a choice of monomial ordering, $|\mathrm{FLAG}_1(K)|$ is an invariant of knots, 
while the classical Alexander $k=1$ ideal is always principal for 
classical knots. We note that setting $s=1-t$ and $ss^{-1}=1$ in each of the 
polynomials in the $\mathrm{FLAG}_1^<$ ideal yields either the Alexander 
polynomial or 0, since this is the ideal which is set to zero when defining 
the Alexander invariant. Since the $\mathrm{FLAG}_k^<$ ideals are in general 
not principal, these Gr\"obner bases in general contain more information than the
usual $k$th Alexander polynomial.

\begin{example}\textup{
We computed the $\mathrm{FLAG}_1^<$ invariants with graded reverse lexicographical 
ordering with respect to the monomial ordering $t<s<t^{-1}<s^{-1}$ for all 
prime classical knots with up to eight crossings using \texttt{Python} code 
available at the first author's website \texttt{www.esotericka.org}. The results 
are collected in the tables below.\footnote{We note that each of the 
$|\mathrm{FLAG}_1(K)|$ values in the table
are in the set $\{4,7\}.$ Whether this is related to the fact that the computation
was done on the Pomona College campus is currently unknown.}
\[\scalebox{0.95}{$
\begin{array}{c|c|l}
K &|\mathrm{FLAG}_1(K)| & \mathrm{FLAG}_1(K) \\ \hline
3_1 & 4 &\{t^2 - t + 1,\  s^{-1} - t,\  t^{-1} + t - 1,\  s + t - 1\} \\ 
4_1 & 4 &\{t^2 - 3t + 1,\  s^{-1} + t - 2,\  ^t{-1} + t - 3,\  s + t - 1\} \\ 
5_1 & 7 &\{s^{-2} - s^{-1} + t^{-1} + t + 1,\  s^{-1}t^{-1} - s^{-1} - t^{-1},\  \\
 & & \quad s^{-1}t - s^{-1} + 1,\  s^{-1} + t^{-2} - t,\  t^{-1}t - 1,\  -s^{-1} - t^{-1} + t^2 + 1,\  s + t - 1\} \\
5_2 & 4 &\{2t^2 - 3t + 2,\  s^{-1} - 2t + 1,\  2t^{-1} + 2t - 3,\  s + t - 1\}\\
6_1 & 4 &\{(t - 2)(2t - 1),\  s^{-1} + 2t - 3,\  2t^{-1} + 2t - 5,\  s + t - 1\}\\ 
6_2 & 7 &\{s^{-2}-s^{-1} - t^{-1} - t + 1,\  s^{-1}t^{-1} - s^{-1} - t^{-1},\  s^{-1}t - s^{-1} + 1,\\
 & & \quad  -s^{-1} + t^{-2} - 2t^{-1} - t + 2,\  t^{-1}t - 1,\  s^{-1} - t^{-1} + t^2 - 2t + 1,\  s + t - 1\}\\ 
6_3 & 7 &\{s^{-2}-s^{-1} + t^{-1} + t - 1,\  s^{-1}t^{-1} - s^{-1} - t^{-1},\  s^{-1}t - s^{-1} + 1,\\
 & & \quad s^{-1} + t^{-2} - 2t^{-1} - t + 2,\  t^{-1}t - 1,\  -s^{-1} - t^{-1} + t^2 - 2t + 3,\  s + t - 1\}\\ 
7_1 & 7 &\{s^{-1}+t^{-3}+t^{-1}-t^2-1,\ -s^{-1}-t^{-2} + t^3 + t,\\
& & \quad   s^{-2} - s^{-1} + t^{-2} + t^{-1} + t^2 + t + 2,\ s^{-1}t^{-1} - s^{-1} - t^{-1},\  s^{-1}t - s^{-1} + 1,\\
& & \quad  t^{-1}t - 1,\  s + t - 1\}\\ 
7_2 & 4 &\{3t^2 - 5t + 3,\  s^{-1} - 3t + 2,\  3t^{-1} + 3t - 5,\  s + t - 1\}\\ 
7_3 & 7 &\{s^{-2} - s^{-1} + 2t^{-1} + 2t + 1,\  s^{-1}t^{-1} - s^{-1} - t^{-1},\  s^{-1}t - s^{-1} + 1,\\
 & & \quad s^{-1} + 2t^{-2} - t^{-1} - 2t + 1,\  t^{-1}t - 1,\  -s^{-1} - 2t^{-1} + 2t^2 - t + 2,\  s + t - 1\}\\ 
7_4 & 4 &\{4t^2 - 7t + 4,\  s^{-1} - 4t + 3,\  4t^{-1} + 4t - 7,\  s + t - 1\}\\ 
7_5 & 7 &\{s^{-2} - s^{-1} + 2t^{-1} + 2t,\  s^{-1}t^{-1} - s^{-1} - t^{-1},\  s^{-1}t - s^{-1} + 1,\\
 & & \quad  s^{-1} + 2t^{-2} - 2t^{-1} - 2t + 2,\  t^{-1}t - 1,\  -s^{-1} - 2t^{-1} + 2t^2 - 2t + 3,\\
& & \quad  s + t - 1\}\\ 
7_6 & 7 &\{s^{-2} - s^{-1} - t^{-1} - t + 3,\  s^{-1}t^{-1} - s^{-1} - t^{-1},\  s^{-1}t - s^{-1} + 1,\\
 & & \quad -s^{-1} + t^{-2} - 4t^{-1} - t + 4,\  t^{-1}t - 1,\  s^{-1} - t^{-1} + t^2 - 4t + 3,\  s + t - 1\}\\ 
7_7 & 7 &\{s^{-2} - s^{-1} + t^{-1} + t - 3,\  s^{-1}t^{-1} - s^{-1} - t^{-1},\  s^{-1}t - s^{-1} + 1,\\
 & & \quad s^{-1} + t^{-2} - 4t^{-1} - t + 4,\  t^{-1}t - 1,\  -s^{-1} - t^{-1} + t^2 - 4t + 5,\  s + t - 1\}\\ 
8_1 & 4 & \{3t^2 - 7t + 3,\  s^{-1} + 3t - 4,\  3t^{-1} + 3t - 7,\  s + t - 1\}\\ 
8_2 & 7 & \{-s^{-1} + t^{-3} - 2t^{-2} + t^{-1} - t^2 + 2t - 1,\  s^{-1} - t^{-2} + 2t^{-1} + t^3 - 2t^2 + t - 2,\\
 & & \quad s^{-2} - s^{-1} - t^{-2} + t^{-1} - t^2 + t,\  s^{-1}t^{-1} - s^{-1} - t^{-1},\  s^{-1}t - s^{-1} + 1,\  \\ 
 & & \quad t^{-1}t - 1,\ s + t - 1\}\\
8_3 & 4 & \{4t^2 - 9t + 4,\  s^{-1} + 4t - 5,\  4t^{-1} + 4t - 9,\  s + t - 1\}\\ 
8_4 & 7 & \{s^{-2} - s^{-1} - 2t^{-1} - 2t + 1,\  s^{-1}t^{-1} - s^{-1} - t^{-1},\  s^{-1}t - s^{-1} + 1,\\
& & \quad   -s^{-1} + 2t^{-2} - 3t^{-1} - 2t + 3,\  t^{-1}t - 1,\  s^{-1} - 2t^{-1} + 2t^2 - 3t + 2,\  s + t - 1\}\\ 
8_5 & 7 & \{-s^{-1} + t^{-3} - 2t^{-2} + 2t^{-1} - t^2 + 2t - 2,\  s^{-1} - t^{-2} + 2t^{-1} + t^3 - 2t^2 + 2t - 3,\\
& & \quad  s^{-2} - s^{-1} - t^{-2} + t^{-1} - t^2 + t - 1,\  s^{-1}t^{-1} - s^{-1} - t^{-1},\\
& & \quad  s^{-1}t - s^{-1} + 1,\  t^{-1}t - 1,\  s + t - 1\}\\ 
8_6 & 7 & \{s^{-2} - s^{-1} - 2t^{-1} - 2t + 2,\  s^{-1}t^{-1} - s^{-1} - t^{-1},\  s^{-1}t - s^{-1} + 1,\\
& & \quad  -s^{-1} + 2t^{-2} - 4t^{-1} - 2t + 4,\  t^{-1}t - 1,\  s^{-1} - 2t^{-1} + 2t^2 - 4t + 3,\  s + t - 1\}\\ 
8_7 & 7 & \{s^{-1} + t^{-3} - 2t^{-2} + 3t^{-1} - t^2 + 2t - 3,\  -s^{-1} - t^{-2} + 2t^{-1} + t^3 - 2t^2 + 3t - 2,\\
&  & \quad  s^{-2} - s^{-1} + t^{-2} - t^{-1} + t^2 - t + 2,\  s^{-1}t^{-1} - s^{-1} - t^{-1},\  s^{-1}t - s^{-1} + 1,\\
&  & \quad  t^{-1}t - 1,\  s + t - 1\}\\ 
8_8 & 7 & \{s^{-2} - s^{-1} + 2t^{-1} + 2t - 2,\  s^{-1}t^{-1} - s^{-1} - t^{-1},\  s^{-1}t - s^{-1} + 1,\\
& & \quad  s^{-1} + 2t^{-2} - 4t^{-1} - 2t + 4,\  t^{-1}t - 1,\  -s^{-1} - 2t^{-1} + 2t^2 - 4t + 5,\  s + t - 1\}\\ 
8_9 & 7  & \{-s^{-1} + t^{-3} - 2t^{-2} + 3t^{-1} - t^2 + 2t - 3,\  s^{-1} - t^{-2} + 2t^{-1} + t^3 - 2t^2 + 3t - 4,\\
& & \quad  s^{-2} - s^{-1} - t^{-2} + t^{-1} - t^2 + t - 2,\  s^{-1}t^{-1} - s^{-1} - t^{-1},\  s^{-1}t - s^{-1} + 1,\\
& & \quad  t^{-1}t - 1,\  s + t - 1\}\\ 
8_{10} & 7 & \{s^{-1} + t^{-3} - 2t^{-2} + 4t^{-1} - t^2 + 2t - 4,\  -s^{-1} - t^{-2} + 2t^{-1} + t^3 - 2t^2 + 4t - 3,\\
& & \quad  s^{-2} - s^{-1} + t^{-2} - t^{-1} + t^2 - t + 3,\  s^{-1}t^{-1} - s^{-1} - t^{-1},\  s^{-1}t - s^{-1} + 1,\\
& & \quad  t^{-1}t - 1,\  s + t - 1\}\\ 
\end{array}$}\]\[\scalebox{0.95}{$\begin{array}{c|c|l} K &|\mathrm{FLAG}_1(K)| & FLAG(K) \\ \hline
8_{11} & 7 & \{s^{-2} - s^{-1} - 2t^{-1} - 2t + 3,\  s^{-1}t^{-1} - s^{-1} - t^{-1},\  s^{-1}t - s^{-1} + 1,\\
& & \quad  -s^{-1} + 2t^{-2} - 5t^{-1} - 2t + 5,\  t^{-1}t - 1,\  s^{-1} - 2t^{-1} + 2t^2 - 5t + 4,\  s + t - 1\}\\ 
8_{12} & 7 & \{s^{-2} - s^{-1} + t^{-1} + t - 5,\  s^{-1}t^{-1} - s^{-1} - t^{-1},\  s^{-1}t - s^{-1} + 1,\\
& & \quad  s^{-1} + t^{-2} - 6t^{-1} - t + 6,\  t^{-1}t - 1,\  -s^{-1} - t^{-1} + t^2 - 6t + 7,\  s + t - 1\}\\ 
8_{13} & 7 &\{s^{-2} - s^{-1} + 2t^{-1} + 2t - 3,\  s^{-1}t^{-1} - s^{-1} - t^{-1},\  s^{-1}t - s^{-1} + 1,\\
& & \quad  s^{-1} + 2t^{-2} - 5t^{-1} - 2t + 5,\  t^{-1}t - 1,\  -s^{-1} - 2t^{-1} + 2t^2 - 5t + 6,\  s + t - 1\}\\ 
8_{14} & 7 &\{s^{-2} - s^{-1} - 2t^{-1} - 2t + 4,\  s^{-1}t^{-1} - s^{-1} - t^{-1},\  s^{-1}t - s^{-1} + 1,\\
& & \quad  -s^{-1} + 2t^{-2} - 6t^{-1} - 2t + 6,\  t^{-1}t - 1,\  s^{-1} - 2t^{-1} + 2t^2 - 6t + 5,\  s + t - 1\}\\ 
8_{15} & 7 &\{s^{-2} - s^{-1} + 3t^{-1} + 3t - 2,\  s^{-1}t^{-1} - s^{-1} - t^{-1},\  s^{-1}t - s^{-1} + 1,\\
& & \quad  s^{-1} + 3t^{-2} - 5t^{-1} - 3t + 5,\  t^{-1}t - 1,\  -s^{-1} - 3t^{-1} + 3t^2 - 5t + 6,\  s + t - 1\}\\ 
8_{16} & 7 & \{s^{-1} + t^{-3} - 3t^{-2} + 5t^{-1} - t^2 + 3t - 5,\  -s^{-1} - t^{-2} + 3t^{-1} + t^3 - 3t^2 + 5t - 4,\\
& & \quad  s^{-2} - s^{-1} + t^{-2} - 2t^{-1} + t^2 - 2t + 3,\  s^{-1}t^{-1} - s^{-1} - t^{-1},\  s^{-1}t - s^{-1} + 1,\\
& & \quad  t^{-1}t - 1,\  s + t - 1\}\\ 
8_{17} & 7 &\{-s^{-1} + t^{-3} - 3t^{-2} + 5t^{-1} - t^2 + 3t - 5,\  s^{-1} - t^{-2} + 3t^{-1} + t^3 - 3t^2 + 5t - 6,\\
& & \quad  s^{-2} - s^{-1} - t^{-2} + 2t^{-1} - t^2 + 2t - 3,\  s^{-1}t^{-1} - s^{-1} - t^{-1},\  s^{-1}t - s^{-1} + 1,\\
& & \quad  t^{-1}t - 1,\  s + t - 1\}\\ 
8_{18} & 7 & \{-s^{-1} + t^{-3} - 4t^{-2} + 6t^{-1} - t^2 + 4t - 6,\  s^{-1} - t^{-2} + 4t^{-1} + t^3 - 4t^2 + 6t - 7,\\
& & \quad  s^{-2} - s^{-1} - t^{-2} + 3t^{-1} - t^2 + 3t - 3,\  s^{-1}t^{-1} - s^{-1} - t^{-1},\  s^{-1}t - s^{-1} + 1,\\
& & \quad  t^{-1}t - 1,\  s + t - 1\}\\ 
8_{19} & 7 & \{s^{-1} + t^{-3} - t^2,\  -s^{-1} - t^{-2} + t^3 + 1,\  s^{-2} - s^{-1} + t^{-2} + t^{-1} + t^2 + t + 1,\\
& & \quad  s^{-1}t^{-1} - s^{-1} - t^{-1},\  s^{-1}t - s^{-1} + 1,\  t^{-1}t - 1,\  s + t - 1\}\\ 
8_{20} & 7 & \{s^{-2} - s^{-1} + t^{-1} + t,\  s^{-1}t^{-1} - s^{-1} - t^{-1},\  s^{-1}t - s^{-1} + 1,\\
& & \quad  s^{-1} + t^{-2} - t^{-1} - t + 1,\  t^{-1}t - 1,\  -s^{-1} - t^{-1} + t^2 - t + 2,\  s + t - 1\}\\  
8_{21} & 7 & \{s^{-2} - s^{-1} - t^{-1} - t + 2,\  s^{-1}t^{-1} - s^{-1} - t^{-1},\  s^{-1}t - s^{-1} + 1,\\
& & \quad   -s^{-1} + t^{-2} - 3t^{-1} - t + 3,\  t^{-1}t - 1,\  s^{-1} - t^{-1} + t^2 - 3t + 2,\  s + t - 1\}
\end{array}$}
\]
}\end{example}

The FLAG ideals are defined for virtual knots just as for classical knots
since each crossing can be considered locally. Thus, we can extend the FLAG invariant
to virtual knots in the usual way by simply ignoring virtual crossings.

\begin{example}\textup{
The virtual knot below, named $4.99$ in the knot atlas \cite{KA}, has trivial virtual 
Alexander polynomial, as does the trefoil $3_1$. However, the two are distinguished by their 
$FLAG_1$ invariants; this shows that the $FLAG_1^<$ invariant (with the same monomial 
ordering as above) is not determined by the virtual Alexander polynomial.
\[\includegraphics{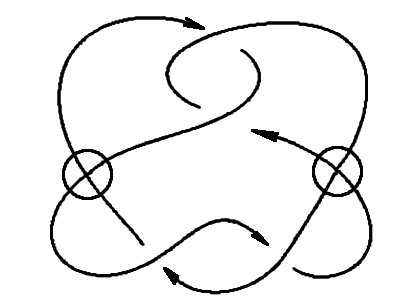} \quad \raisebox{0.5in}{$\mathrm{FLAG_1}(4.99)
=\{s^{-1}-2,t^{-1}-2,2s-1,2t-1\}$}\]
We remark that this virtual knot has classical Alexander polynomial $2t-1$, which is not
symmetric, unlike the case for all classical knots.
}\end{example}

\section{\large \textbf{Questions}}\label{QFR}

In this section we collect a few questions for future research.

What other quotients of the fundamental quandle yield interesting finite quandles?
What is the relationship between quotients of the fundamental quandle, a variety
of functorial invariant, and the homomorphism-based invariants such as the quandle
counting invariant and its enhancements? What does the cardinality of the 
$\mathrm{FLAG}_k^<$ invariant tell us about a knot?

\bigskip

\noindent\textsc{Department of Mathematics \\ 
Claremont McKenna College \\
850 Columbia Ave. \\
Claremont, CA 91711}

\end{document}